\documentclass[11pt]{article}

\usepackage{amsmath}
\usepackage{amssymb}
\usepackage{eucal}

\begin{document}

\baselineskip 16pt

\title{On $\sigma$-quasinormal  subgroups of finite groups
 \thanks{Research is supported by an NNSF grant of China (Grant No.
11401264) and a TAPP of Jiangsu Higher Education Institutions (PPZY 2015A013)}}

\author{ Bin Hu\\\
{\small School of Mathematics and Statistics, Jiangsu Normal University,}\\ {\small Xuzhou, 221116, P.R. China}\\
{\small E-mail: hubin118@126.com}\\ \\
Jianhong Huang  \thanks{Corresponding author}\\
{\small School of Mathematics and Statistics, Jiangsu Normal University,}\\
{\small Xuzhou 221116, P. R. China}\\
{\small E-mail: jhh320@126.com}\\ \\
{ Alexander  N. Skiba   }\\
{\small Department of Mathematics  and Technologies of Programming,  Francisk Skorina Gomel State University,}\\
{\small Gomel 246019, Belarus}\\
{\small E-mail: alexander.skiba49@gmail.com}}

\date{}
\maketitle

\begin{abstract}    Let  $G$ be a finite group and
 $\sigma =\{\sigma_{i} | i\in I\}$  some
partition of the set  of all primes $\Bbb{P}$, that is,  $\sigma =\{\sigma_{i} |
 i\in I \}$, where   $\Bbb{P}=\bigcup_{i\in I} \sigma_{i}$
 and $\sigma_{i}\cap
\sigma_{j}= \emptyset  $ for all $i\ne j$.   We say that $G$  is 
\emph{$\sigma$-primary} if $G$ is a  $\sigma _{i}$-group for some $i$.
A subgroup $A$ of $G$ is said to be: \emph{${\sigma}$-subnormal}
  in $G$
 if   there is a subgroup chain  $A=A_{0} \leq A_{1} \leq \cdots \leq
A_{n}=G$  such that  either $A_{i-1}\trianglelefteq A_{i}$ or
$A_{i}/(A_{i-1})_{A_{i}}$ is   $\sigma$-primary
  for all $i=1, \ldots , n$;
 \emph{modular} in $G$  if the following conditions hold:
(i) $\langle X, A \cap Z \rangle=\langle X, A \rangle \cap Z$ for all $X \leq G, Z \leq
 G$ such that
$X \leq Z$, and
(ii) $\langle A, Y \cap Z \rangle=\langle A, Y \rangle \cap Z$ for all $Y \leq G, Z \leq
 G$ such that
$A \leq Z$.

In this paper, a subgroup $A$ of $G$ is   called     \emph{$\sigma$-quasinormal in
$G$}
 if  $L$ is modular 
 and    ${\sigma}$-subnormal in $G$.

We study  $\sigma$-quasinormal subgroups of $G$. In
particular, we prove that  if a subgroup $H$ of 
 $G$ is $\sigma$-quasinormal
in $G$, then for every chief factor $H/K$ of $G$ between $H^{G}$ and
$H_{G}$ the semidirect product $(H/K)\rtimes (G/C_{G}(H/K))$ is
$\sigma$-primary.

\end{abstract}


\let\thefootnote\thefootnoteorig

\section{Introduction}

Throughout this paper, all groups are finite and $G$ always denotes
a finite group. Moreover,  $\mathbb{P}$ is the set of all  primes,
  $\pi= \{p_{1}, \ldots , p_{n}\} \subseteq  \Bbb{P}$ and  $\pi' =  \Bbb{P} \setminus \pi$. If
 $n$ is an integer, the symbol $\pi (n)$ denotes
 the set of all primes dividing $n$; as usual,  $\pi (G)=\pi (|G|)$, the set of all
  primes dividing the order of $G$.

A  subgroup $A$ of $G$ is said to be  
 \emph{modular in $G$} \cite{1-3} if it is 
a  modular element
(in the sense of
 Kurosh \cite[p. 43]{Schm})  of the
 lattice of all subgroups of $G$, that is, the following conditions hold:

(i) $\langle X,M \cap Z \rangle=\langle X, M \rangle \cap Z$ for all $X \leq G, Z \leq
 G$ such that
$X \leq Z$, and

(ii) $\langle M, Y \cap Z \rangle=\langle M, Y \rangle \cap Z$ for all $Y \leq G, Z \leq
 G$ such that
$M \leq Z$.

In what follows, $\sigma$  is some partition of
$\Bbb{P}$, that is,  $\sigma =\{\sigma_{i} |
 i\in I \}$, where   $\Bbb{P}=\bigcup_{i\in I} \sigma_{i}$
 and $\sigma_{i}\cap
\sigma_{j}= \emptyset  $ for all $i\ne j$.  
By the analogy with the notation   $\pi (n)$, we write  $\sigma (n)$ to denote 
the set  $\{\sigma_{i} |\sigma_{i}\cap \pi (n)\ne 
 \emptyset  \}$;   $\sigma (G)=\sigma (|G|)$.

 The group  $G$ is said to be:
\emph{$\sigma$-primary} \cite{1} if  $G$ is a $\sigma_{i}$-group for some 
$i$; \emph{$\sigma$-decomposable} (Shemetkov \cite{Shem})  or
\emph{$\sigma$-nilpotent}  (Guo and Skiba \cite{33}) if $G=G_{1}\times \dots \times G_{n}$
for some $\sigma$-primary groups $G_{1}, \ldots, G_{n}$.  We use $\mathfrak{N}_{\sigma}$ to denote
 the class of all  $\sigma$-nilpotent groups.

We say, following \cite{1}, that the  subgroup $A$ of $G$ is
 \emph{${\sigma}$-subnormal
  in $G$}  if  it  is \emph{$\mathfrak{N}_{\sigma}$-subnormal in $G$} in the
 sense of Kegel  \cite{KegSubn}, that is, there is a subgroup chain  $$A=A_{0} \leq A_{1} \leq \cdots \leq
A_{n}=G$$  such that  either $A_{i-1} \trianglelefteq A_{i}$ or
$A_{i}/(A_{i-1})_{A_{i}}$ is  ${\sigma}$-nilpotent 
  for all $i=1, \ldots , n$.

 A  subgroup $A$ of  $G$ is said to be  \emph{ quasinormal } \cite{5} or \emph{ permutable}
 \cite{6, prod} in $G$    if  $A$ permutes with every subgroup $L$ of $G$, that is, $AL=LA$.

The quasinormal subgroups have many interesting properties.
For instance, if $A$ is quasinormal in $G$, then:  {\sl $A$ is subnormal in $G$}
 (Ore \cite{5}),  {\sl $A/A_{G}$ is nilpotent} (Ito and Szep \cite{It}) and,
 in general,  \emph{$A/A_{G}$ is not necessarily abelian}
(Thomson  \cite{Th}). Every quasinormal subgroup $A$ of $G$ is modular in  
$G$ \cite{1-3}.
 Moreover, the following properties of  quasinormal subgroups are well-known.

{\bf Theorem A} (See Theorem 5.1.1  in \cite{Schm}).   {\sl A subgroup $A$ of $G$ is
quasinormal in $G$ if and only if $A$ is  modular and subnormal in
$G$. }

{\bf Theorem B.} {\sl  If $A$ is
 a quasinormal subgroup of $G$, then:}

(i)   {\sl $A^{G}/A_{G}$ is nilpotent}  (This follows from the above-mentioned results
 in \cite{5, It}), {\sl      and }

(ii) {\sl   Every chief factor $H/K$ of $G$ between $A^{G}$ and
$A_{G}$ is central in $G$, that is, $C_{G}(H/K)=G$ }   (Maier and Schmid
\cite{MaierS}).

Since every subnormal subgroup of $G$ is $\sigma$-subnormal  in $G$,
 Theorems A and B make natural to ask: {\sl What we can say about
 the quotient $A^{G}/A_{G}$ provided
the subgroup $A$ of $G$ is  $\sigma$-qusinormal in the sense of the
following definition?}

{\bf Definition 1.1.}  Let $A$ be a subgroup of $G$. Then we say that $A$
is   \emph{$\sigma$-qusinormal} in $G$  if $A$ is    modular and
  $\sigma$-subnormal in $G$.

In this note we give the answer to this question. But before continuing,
consider the following

{\bf Example 1.2.}    Let $p >  q,  r, t$ be  distinct  primes, where
 $t$ divides $r-1$. Let     $Q$ be   a simple
  ${\mathbb F}_{q}C_{p}$-module  which is     faithful  for $C_{p}$, let 
$C_{r}\rtimes C_{t}$   be a non-abelian group of order $rt$, and let $A=C_{t}$.
 Finally, let   $G=(Q\rtimes C_{p})\times (C_{r}\rtimes C_{t})$ and   $B$ be
  a subgroup of order $q$ in $Q$. Then $ B < Q$ since $p  > q$.
It is not difficult  to show that  $A$ is modular in $G$ (see  \cite[Lemma 5.1.8]{Schm}).
 On the other hand,
  $A$ is   $\sigma$-subnormal  in $G$, where
 $\sigma =\{\{q, r, t\}, \{q, r, t\}'\}$.  Hence
  $A$ is $\sigma$-quasinormal   in $G$. It is clear also that $A$ is not subnormal in $G$, so
$A$ is not quasinormal in $G$. Finally, note that $B$ is not modular in
$G$ by Lemma 2.2 below.

A  chief factor $H/K$ of $G$ is said to be \emph{$\sigma$-central} in $G$  \cite{commun} if
the semidirect product  $(H/K)\rtimes (G/C_{G}(H/K))$ is $\sigma$-primary. Note that
$G$ is $\sigma$-nilpotent if and only if every chief factor of $G$ is $\sigma$-central in $G$
 \cite{1}.  A subgroup $A$ of $G$ is said to be: \emph{$\sigma$-seminormal  
in $G$} (J.C. Beidleman) if $x\in N_{G}(A)$ for all $x\in G$ such that 
 $\sigma (|x|)\cap \sigma (A)=\emptyset$; \emph{seminormal  
in $G$} if $x\in N_{G}(A)$ for all $x\in G$ such that 
 $\pi (|x|)\cap \pi (A)=\emptyset$.

Our main goal here is to prove the following

{\bf Theorem C.} {\sl  Let $A$ be a $\sigma$-quasinormal subgroup of   $G$.
 Then the following
statements hold:}

(i) {\sl If $G$
possesses a  Hall $\sigma_{i}$-subgroup,
 then $A$ permutes with each Hall $\sigma_{i}$-subgroup  of  $G$. }

(ii) {\sl The quotients
$A^{G}/A_{G}$ and $G/C_{G}(A^{G}/A_{G}) $ are $\sigma$-nilpotent, and }

(iii) {\sl Every chief factor of $G$ between $A^{G}$ and $A_{G}$ is $\sigma$-central in $G$.  }

(iv)  {\sl For every $i$ such that $\sigma _{i} \in \sigma 
(G/C_{G}(A^{G}/A_{G}))$  we have 
 $\sigma _{i} \in  \sigma (A^{G}/A_{G}).$
}

(v) {\sl $A$  is $\sigma$-seminormal in $G$.}

 The subgroup $A$ of $G$ is subnormal in $G$ if and only if it is
$\sigma$-subnormal in $G$, where  $\sigma=\sigma ^{1} =\{\{2\}, \{3\},
\ldots \}$ (we use here the terminology in \cite{alg12}).  It is clear also
 that $G$ is nilpotent if and only if  $G$ is $
\sigma ^{1}$-nilpotent, and a  chief factor $H/K$ of $G$ is central in $G$ if
and only if  $H/K$ is $ \sigma ^{1}$-central in $G$. Therefore Theorem B
is a special case of Theorem C, when  $\sigma =  \sigma ^{1}$.

 In the other classical case when  $\sigma =\sigma ^{\pi}=\{\pi, \pi'\}$:
 $G$  is  $\sigma ^{\pi}$-nilpotent 
 if and only if $G$ is  \emph{ $\pi$-decomposable}, that is,
 $G=O_{\pi}(G)\times O_{\pi'}(G)$); a subgroup  $A$ of 
 $G$ is $\sigma ^{\pi}$-subnormal in $G$ if and only if there is a subgroup
 chain  $$A=A_{0} \leq A_{1} \leq \cdots \leq
A_{n}=G$$  such that either $A_{i-1} \trianglelefteq A_{i}$ or 
 $A_{i}/(A_{i-1})_{A_{i}}$ is a  $\pi _{0}$-group, where   $\pi _{0}\in  \{\pi,
\pi'\}$,    for all $i=1, \ldots , n$.
 Thus,    in this case we get from Theorem C the following

{\bf Corollary 1.3.} {\sl  Suppose that  $A$ is modular subgroup of   $G$ and there is
 a subgroup chain  $$A=A_{0} \leq A_{1} \leq \cdots \leq
A_{n}=G$$  such that  either $A_{i-1} \trianglelefteq A_{i}$ or
$A_{i}/(A_{i-1})_{A_{i}}$ is a  $\pi _{0}$-group, where   $\pi _{0}\in  \{\pi,
\pi'\}$,
  for all $i=1, \ldots , n$.
 Then the following statements hold:}

(i) {\sl If $G$
possesses a  Hall $\pi _{0}$-subgroup, where   $\pi _{0}\in  \{\pi,
\pi'\}$,
 then $A$ permutes with each Hall $\pi _{0}$-subgroup  of  $G$. }

(ii) {\sl The quotients
$A^{G}/A_{G}$ and $G/C_{G}(A^{G}/A_{G}) $ are $\pi$-decomposable, and }

(iii) {\sl For every  chief
 factor $H/K$ of $G$ between $A^{G}$ and $A_{G}$
the semidirect product $(H/K)\rtimes (G/C_{G}(H/K))$ is either a $\pi$-group or a $\pi'$-group.  }

In fact, in the theory of $\pi$-soluble groups ($\pi= \{p_{1}, \ldots , p_{n}\}$)
 we deal with the  partition 
$\sigma =\sigma ^{1\pi }=\{\{p_{1}\}, \ldots , \{p_{n}\}, \pi'\}$ of $\Bbb{P}$. 
Note that  $G$ is $\sigma ^{1\pi }$-nilpotent
 if and only if  $G$ is    \emph{$\pi$-special} \cite{Cun2},
 that is, $G=O_{p_{1}}(G)\times \cdots \times
 O_{p_{n}}(G)\times O_{\pi'}(G)$.  A subgroup  $A$ of $G$ is $\sigma ^{1\pi }$-subnormal
 in $G$ if and only if it is
$\frak{F}$-subnormal in $G$ in the sense of Kegel \cite{KegSubn}, where  
$\frak{F}$ is  the class of all $\pi'$-groups.   Therefore, in this case
 we get from Theorem C the following

{\bf Corollary 1.4.} {\sl  Suppose that  $A$ is a modular subgroup of   $G$ and there is
 a subgroup chain  $$A=A_{0} \leq A_{1} \leq \cdots \leq
A_{n}=G$$  such that  either $A_{i-1} \trianglelefteq A_{i}$ or
$A_{i}/(A_{i-1})_{A_{i}}$ is a  $\pi'$-group. Then the following
statements hold:}

(i) {\sl   $A$ permutes with every Sylow $p$-subgroup of $G$ for all $p\in
\pi$, and if $G$
possesses a  Hall $\pi'$-subgroup, 
 then $A$ permutes with each  Hall $\pi'$-subgroup  of  $G$. }

(ii) {\sl The quotients
$A^{G}/A_{G}$ and $G/C_{G}(A^{G}/A_{G}) $ are $\pi$-special, and }

(iii) {\sl For every   non-central   chief
 factor $H/K$ of $G$ between $A^{G}$ and $A_{G}$
the semidirect product $(H/K)\rtimes (G/C_{G}(H/K))$ is a $\pi'$-group.  }

\section{Preliminaries}

If $G=A\rtimes \langle t \rangle$ is non-abelian  with an elementary abelian
$p$-group $A$ and an element $t$ of
 prime order $q\ne p$ induces a non-trivial power
 automorphism  on $A$, then we say  $G$ is a \emph{$P$-group of type  $(p, q)$}
(see \cite[p. 49]{Schm}).

{\bf Lemma 2.1}  (See Lemma 2.2.2(d)   in \cite{Schm}).
 {\sl If  $G=A\rtimes \langle t \rangle$ is a $P$-group of type  $(p, q)$, then
 $\langle t \rangle^{G}=G$.}

The following remarkable result of R. Schmidt plays a key role in the 
proof of Theorem C. 

{\bf Lemma 2.2}  (See Theorems 5.1.14 in \cite{Schm}).
 {\sl Let $M$  be a modular subgroup of $G$ with $M_{G}=1$.   Then 
 $ G=S_{1}\times \cdots \times S_{r} \times K,$
 where  $0\leq r\in \mathbb{Z}$ and for all $i, j \in \{1, \ldots , r\}$,
}

(a) {\sl $S_{i}$ is a non-abelian $P$-group,}

(b) {\sl  $(|S_{i}|, |S_{j}|)=1=(|S_{i}|, |K|)$ for all $i\ne j$, }

(c) {\sl $M=Q_{1}\times \cdots \times Q_{r}\times (M\cap K)$ and $Q_{i}$
 is a non-normal Sylow subgroup of $S_{i}$, }

(d) {\sl $M\cap K$ is quasinormal in $G$. }

{\bf Lemma 2.3} (See Lemma 2.6 in \cite{1}). {\sl Let $A$, $B$  and
$N$ be subgroups of $G$, where $A$ is $\sigma$-subnormal  and $N$ is
normal in $G$. Then:}

(1) {\sl $A \cap B$   is $\sigma$-subnormal   in  $B$.}

(2) {\sl   $AN/N$  is $\sigma$-subnormal in $G/N$.}

(3) {\sl  $A\cap H$  is a 
 Hall $\sigma _{i}$-subgroup of $A$  for every  Hall $\sigma _{i}$-subgroup $H$ of $G$. }

The following lemma is a corollary of  Lemma 2.3 and  
general properties  of modular subgroups \cite[p. 201]{Schm}.

{\bf Lemma 2.4.}  {\sl Let $A$, $B$  and
$N$ be subgroups of $G$, where $A$ is $\sigma$-quasinormal and $N$ is
normal in $G$.}

(1) {\sl If $A \leq  B$, then  $A$ is $\sigma$-quasinormal  in   $B$.}
 
(2) {\sl  $AN/N$
 is $\sigma$-quasinormal  in  $G/N$}.

A normal subgroup $E$ of $G$ is said to be \emph{${\sigma}$-hypercentral}  (in
$G$) if either $E=1$ or   every chief factor of $G$ below $E$ is $\sigma$-central.
We use  $Z_{\sigma}(G)$ to denote the  \emph{${\sigma}$-hypercentre } of $G$ \cite{commun},
 that is,
 the product of all normal
${\sigma}$-hypercentral  subgroups of $G$.

{\bf Lemma   2.5. }  {\sl    Every  chief factor
 of $G$ below $Z_{\sigma}(G)$ is $\sigma$-central in $G$. }

{\bf Proof. }   It is enough to consider the case when
$Z=A_{1}A_{2}$, where
$A_{1}$ and $A_{2}$ are normal ${\sigma}$-hypercentral subgroups of
 $G$. Moreover, in view of the  Jordan-H\"{o}lder theorem for the chief series,
it is enough to show that if   $A_{1}\leq K < H \leq A_{1}A_{2}$, then $H/K$
   is $\sigma$-central. But in this case we
have  $H=A_{1}(H\cap A_{2})$, where  $H\cap A_{2}\nleq  K$ and so 
from  
the $G$-isomorphism $(H\cap A_{2})/(K\cap A_{2})\simeq  (H\cap 
A_{2})K/K=H/K$ we get that $C_{G}(H/K)=C_{G}((H\cap A_{2})/(K\cap A_{2}))$
  and hence   $H/K$   is $\sigma$-central in $G$.    The lemma is proved.

{\bf Lemma 2.6.} {\sl Let   $N$ be   a normal $\sigma _{i}$-subgroup of   $G$.
 Then   $N\leq  Z_{\sigma}(G)$ if and only if
$O^{\sigma _{i}}(G)\leq C_{G}(N)$.  }

{\bf Proof. } If   $O^{\sigma _{i}}(G)\leq C_{G}(N)$, then
for every chief factor  $H/K$ of $G$ below $N$ both groups $H/K$ and $G/C_{G}(H/K)$ are
${\sigma _{i}}$-group since $G/O^{\sigma _{i}}(G)$ is a ${\sigma _{i}}$-group. Hence
$(H/K)\rtimes (G/C_{G}(H/K))$ is $\sigma$-primary.  Thus   $N\leq  Z_{\sigma}(G)$.

Now assume that $N\leq  Z_{\sigma}(G)$.  Let $1= Z_{0} < Z_{1}
< \cdots <  Z_{t} = N$ be a chief
 series of $G$  below $N$ and  $C_{i}=
C_{G}(Z_{i}/Z_{i-1})$. Let $C= C_{1} \cap \cdots \cap C_{t}$.  Then $G/C$
 is a ${\sigma _{i}}$-group. On the other hand,
 $C/C_{G}(N)\simeq A\leq \text{Aut}(N)$ stabilizes
the series    $1= Z_{0} < Z_{1}
< \cdots <  Z_{t} = N$, so   $C/C_{G}(N)$ is a $\pi (N)$-group by
  \cite[Ch. A, Corollary 12.4(a)]{DH}.   Hence $G/C_{G}(N)$ is  a ${\sigma _{i}}$-group and
 so  $O^{\sigma _{i}}(G)\leq C_{G}(N)$.
 The lemma is proved.

\section{Proof of Theorem C}

   Suppose that this theorem    is false and let $G$
be a counterexample of minimal order.  Then $1  <  A  <  G$. We can assume 
without loss of generality that $\sigma (A)=\{\sigma _{1}, \ldots , \sigma _{m}\}$.

(1) {\sl Statement   (i) holds for $G$.}

   Suppose that this assertion is  false.
 Then for some  Hall $\sigma_{i}$-subgroup $V$
 of   $G$ we have $AV\ne VA$.    It  is  clear that     $V$ is a  Hall $\sigma_{i}$-subgroup $V$
 of $\langle A, V \rangle$.  On the other hand, $A$ is
$\sigma$-quasinormal in  $\langle A, V \rangle$ by Lemma 2.4(1). Therefore
in the case when $\langle A, V \rangle  < G$, we have $AV=VA$ by the
choice of $G$. Thus  $\langle A, V \rangle =G$.

Since  $A$ is $\sigma$-subnormal in $G$,
 there is a subgroup chain  $A=A_{0} \leq A_{1} \leq \cdots \leq
A_{n}=G$  such that  either $A_{i-1}\trianglelefteq A_{i}$ or $A_{i}/(A_{i-1})_{A_{i}}$
 is  ${\sigma}$-primary for all $i=1, \ldots , n$.  We can assume without
loss of generality that $M=A_{n-1} < G$. Then $A$ permutes with every
  Hall $\sigma_ {i}$-subgroup of $M$ by the choice of $G$.

 Moreover, the modularity of $A$ in $G$
implies that $$M=M\cap \langle A,   V \rangle=\langle A, M\cap
V\rangle .$$ On the other hand, by Lemma 2.3(3),  $M\cap V$ is a Hall $\sigma
_{i}$-subgroup of $M$. Hence  $M= A(M\cap V)=(M\cap V)A$ by the choice 
of $G$.

 If  $V\leq  M_{G}$, then
$A  (M\cap V)=A V=VA$ and so $V\nleq
M_{G}$.    
 Now note that  $VM=MV$. Indeed, if $M$ is normal in $G$,
it is clear. Otherwise, $G/M_{G}$ is $\sigma$-primary and so $G=MV=VM$ since
$V\nleq
M_{G}$ and $V$ is a Hall $\sigma_{i}$-subgroup
 of $G$.
  Therefore
 $$VA=V(M\cap V)A=VM=MV
=A(M\cap V)V= AV.$$ This contradiction completes the
proof of (1).

(2) $A_{G}=1$.

Suppose that $A_{G}\ne 1$ and let
$R$   be a minimal normal subgroup of $G$  contained in $A_{G}$. 
 Then $A/R$ is $\sigma$-quasinormal in $G/R$ by  Lemma 2.4(2), so  the
hypothesis holds for $(G/R, A/R)$. Therefore the choice of $G$ implies that
 Statements (ii)--(v)
 hold for $(G/R, A/R)$.   Hence
$$(A/R)^{G/R}/(A/R)_{G/R}=(A^{G}/R)/(A_{G}/R)\simeq A^{G}/A_{G}$$ and
$$(G/R)/C_{G/R}((A/R)^{G/R}/(A/R)_{G/R})=
 (G/R)/(C_{G}(A^{G}/A_{G})/R)\simeq G/C_{G}(A^{G}/A_{G})$$
are   $ \sigma $-nilpotent,
  so Statement  (ii) holds for $G$.

Now let $T/L$ be any   chief factor of $G$ between $A^{G}$ and  $  A_{G}$.
Then $(T/R)/(L/R)$ is a  chief factor of $G/R$ between  $(A/R)^{G/R}$
and  $ (A/R)_{G/R}$  and  so 
   $(T/R)/(L/R)$  is  $\sigma$-central in $G/R$,
that is,
$$((T/R)/(L/R))\rtimes
 ((G/R)/C_{(G/R)}((T/R)/(L/R)))$$ is $\sigma$-primary.
Since the factors  $(T/R)/(L/R)$ and $T/L$ are $G$-isomorphic, it follows
 that $(T/L)\rtimes (G/C_{G}(T/L))$  is $\sigma$-primary too. Hence
$T/L$ is $\sigma$-central in $G$. Thus Statement   (iii) holds
 for $G$.

If  $i$ such that  $$\sigma _{i} \cap \pi (G/C_{G}(A^{G}/A_{G}))  =\sigma _{i} \cap \pi
 ((G/R)/C_{G/R}((A/R)^{G/R}/(A/R)_{G/R}))\ne
 \emptyset,$$  then 
$$\sigma _{i} \cap
 \pi (A^{G}/A_{G}) =\sigma _{i} \cap \pi ((A/R)^{G/R}/(A/R)_{G/R})\ne
\emptyset$$  
 and so Statement (iv) holds for $G$ too.

Finally, if    $x\in G$ is such that $\sigma (A)\cap \sigma (|x|)=\emptyset$,
 then $\sigma (A/R)\cap \sigma (|xR|)=\emptyset$,  so $xR\in N_{G/R}(A/R)=N_{G}(A)/R$ and hence 
Statement (v) holds for $G$.  
 Therefore, in view of Claim (1), the
conclusion of the theorem holds for $G$, which contradicts the choice of
$G$. Hence $A_{G}=1$.

(3) {\sl If $A$ is a $\sigma _{i}$-group, then $A\leq O_{\sigma _{i}}(G)$. 
}

 It is enough to show  that if $A$ is any $\sigma $-subnormal $\sigma _{i}$-subgroup of 
 $G$, then 
$A\leq O_{\sigma _{i}}(G)$.  
Assume that this  is false and let $G$ be a counterexample 
of  minimal order. Then $1 <  A < G$. Let $D=O_{\sigma _{i}}(G)$,  $R$ be 
  a minimal normal subgroup of $G$ and $O/R=O_{\sigma_{i}}(G/R)$.   Then 
the choice of $G$ and Lemma 2.4(ii)  imply that   $AR/R\leq O/R$. 
 Therefore   $R\nleq D $, so $D=1$  and  $A\cap R < R$. It is clear also that 
 $O^{\sigma _{i}}(R)=R$.
 Suppose that $L= A\cap  R\ne 1$.  Lemma 2.3(2)  implies that $L$ 
is $\sigma $-subnormal in $R$.   If $R < G$, the choice of $G$ implies that $L\leq 
O_{\sigma _{i}}(R)\leq   D$ since $O_{\sigma _{i}}(R)$ is a characteristic subgroup of $R$.
   But then $D\ne 1$,
 a contradiction. Hence $R=G$ is a simple group, which is  also 
 impossible since $1 <  A < G$.  Therefore $R\cap 
A=1$.  If  $O < G$, the choice of $G$ implies that 
$A\leq 
O_{\sigma _{i}}(O)\leq   D=1$. Therefore 
  $G/R=O/R$ is a $\sigma _{i}$-group. 
Hence $R$ is a unique minimal normal subgroup of $G$.
 It is clear also  that
 $R\nleq\Phi (G)$,  so $C_{G}(R)\leq R$ by \cite[Ch. A, 15.2]{DH}.

Now we show that  $G=RA$. Indeed, if $RA  < G$, then the choice of $G$ and 
Lemma 2.3(1) imply that  $A\leq O_{\sigma 
_{i}}(RA)$ and so $A= O_{\sigma 
_{i}}(RA)$ since $O_{\sigma _{i}}(R)=1$, which implies that $RA=R\times A$. 
But then $A\leq C_{G}(R)\leq R$ and so $A=1$ since $A\cap R=1$. This 
contradiction shows that $G=RA$.  
 
Since   $A$  is $\sigma 
$-subnormal in  $G$, there is a subgroup $M$ such that $A\leq M  < G$ and either 
$M \trianglelefteq G$ or
$G/M_{G}$ is  ${\sigma}$-primary.
 Since $R$ is a unique minimal normal subgroup of $G$ and $A\leq M  < G=RA$,  $R\nleq M$ and
 hence  $G/M_{G}=G/1$ is a ${\sigma}_{i}$-group.  Therefore  $A\leq O_{\sigma _{i}}(G)=G$.
 This contradiction completes the proof of (3).

(4) {\sl   $A\leq O_{\sigma_{1}}(G)
 \times \cdots \times O_{\sigma _{m}}(G)$.  Hence  $A^{G}=O_{\sigma_{{1}}}(A^{G})
 \times \cdots \times O_{\sigma _{{m}}}(A^{G})$.}

 Claim (2)  and Lemma 2.2(c)(d) imply that  $A=A_{1}\times \cdots \times A_{m}$,
 where  $A_{i}$
is a Hall $\sigma _{i}$-subgroup of $A$ for all $i=1, \ldots , m$.
On the other hand, since $A$ is
$\sigma $-subnormal in $G$, we have  $A_{i}\leq O_{\sigma _{i}}(G)$ by Claim (3). 
Hence we have (4).

(5) {\sl   Statement (iii) holds for $G$.}

Let  $T/L$ be any  chief factor of $
G$ below $A^{G}$.
Suppose  that $T/L$ is not $\sigma$-central in $G$.
 Theorem B(ii) implies then that $A$ is not
quasinormal in $G$, so  in view of
 Lemma 2.2, we have
  $G=S_{1}\times \cdots \times S_{r} \times K,$
 where for all $i, j \in \{1, \ldots , r\}$ the following hold:

(a) {\sl $S_{i}$ is a non-abelian $P$-group,}

(b) {\sl  $(|S_{i}|, |S_{j}|)=1=(|S_{i}|, |K|)$ for  $i\ne j$, }

(c) {\sl $A=Q_{1}\times \cdots \times Q_{r}\times (A\cap K)$ and $Q_{i}$
 is a non-normal Sylow subgroup of $S_{i}$, and }

(d) {\sl $A\cap K$ is quasinormal in $G$. }

Hence, in view of Claim (3),  $$A^{G}=Q_{1}^{G}\times \cdots \times Q_{r}^{G}\times
(A\cap K)^{G}=O_{\sigma_{{1}}}(A^{G})
 \times \cdots \times O_{\sigma _{{m}}}(A^{G}), $$ where  $(A\cap K)^{G} \leq Z_{\infty}(G)\leq Z_{\sigma}(G)$ by Theorem 
B(ii) since  $(A\cap K)_{G}\leq  A_{G}=1$ by Claim (2). 
Therefore, in view of the  Jordan-G\"{o}lder theorem
for the chief series, we can assume  without loss of generality that 
$T\leq S_{k}$ for some $k$.

Now  note  that for  all $i, j$ we have either $S_{i}\leq O_{\sigma_{j}}(A^{G})$ 
  or  $S_{i}\cap O_{\sigma_{j}}(A^{G})=1$.
 Indeed, assume that   $S_{i}\cap O_{\sigma_{j}}(A^{G})\ne 1$. It is clear that for some
$t$ we have $Q_{i}\leq O_{\sigma_{t}}(A^{G})$.  Then $ Q_{i}^{G}=S_{i}
\leq O_{\sigma_{{t}}}(A^{G})$ by Lemma
2.1. Hence $j=t$
 since  $O_{\sigma_{{j}}}(A^{G})\cap O_{\sigma_{{t}}}(A^{G})=1$ for $j\ne t$. Therefore all
$S_{i}$ are $\sigma$-primary.
 Moreover, if $S_{i}$ is a $\sigma _{l}$-group, then
    $G/C_{G}(S_{i})$ is a $\sigma _{l}$-group since  $G=S_{1}\times \cdots
\times S_{r} \times K.$
 Therefore $S_{k}\leq Z_{\sigma}(G)$ by Lemma 2.6 and  so  $T/L$ is
$\sigma$-central in $G$ by Lemma 2.5, a contradiction. 
  Hence  Statement (iii) holds for $G$.

(6) {\sl   Statements (ii) and (iv) hold  for $G$.}

 From Claim   (3) we know that $A^{G} =O_{\sigma_{1}}(A^{G})
 \times \cdots \times O_{\sigma _{m}}(A^{G}).$
 Then $$C_{G}(A^{G})= C_{G}(O_{\sigma_{1}}(A^{G})
) \cap \cdots  \cap C_{G}(O_{\sigma_{m}}(A^{G})
).$$
From Claims (2),  (4)  and Lemma 2.6 we know that $G/C_{G}(O_{\sigma_{i}}(A^{G}))$ is
 a $\sigma _{i}$-group for all $i=1, \dots, m$.
 Therefore, in view of \cite[Ch. I, 9.6]{hupp},
 $$G/C_{G}(A^{G})=G/(C_{G}(O_{\sigma_{1}}(A^{G})) \cap \cdots  \cap
 C_{G}(O_{\sigma_{m}}(A^{G})))$$$$\simeq V\leq (G/C_{G}(O_{\sigma_{1}}(A^{G}))) \times
 \cdots \times (G/C_{G}(O_{\sigma_{m}}(A^{G}))) $$ is  
$\sigma$-nilpotent, and for every $i$ such that $\sigma _{i} \in \sigma 
(G/C_{G}(A^{G}))$  we have 
 $\sigma _{i} \in  \sigma (A^{G}).$ Hence Statements (ii) and  (iv) hold for $G$.

(7) {\sl   Statement (v)  holds  for $G$.}

 Suppose that $x\in G$ is such that $\sigma (A)\cap \sigma 
(|x|)=\emptyset$.  Then  the modularity of $A$ and Claim (4)
 imply that $A=O_{\sigma_{1}}(A^{G})
 \times \cdots \times O_{\sigma _{m}}(A^{G})\cap \langle A, \langle x \rangle   \rangle $
 is normal in 
$\langle A, \langle x \rangle   \rangle $, so $x \in N_{G}(A)$.  Hence we have (7).

 From Claims (1), (5)--(7) it follows that the 
conclusion of the theorem holds for $G$, which contradicts the choice of
$G$. The theorem is proved.

\end{document}